\documentclass{article}

\usepackage{graphicx}
\usepackage{amssymb}
\usepackage{epstopdf}
\usepackage{verbatim}

\usepackage{enumerate}
\usepackage{amsmath}

\usepackage{amsthm}

\newtheorem{lemma}{Lemma}
\newtheorem{theorem}{Theorem}

\newtheorem{definition}{Definition}

\newtheorem{claim}{Claim}
\usepackage{parskip}
\usepackage{cmbright}

\begin{document}

\title{Reconstructing surface triangulations by their intersection matrices}
%\tnotetext[t1]{Research supported by PAPIIT-M\'{e}xico IN112511 and CONACyT 166951}

\author{Jorge Arocha\thanks{Instituto de Matem\'{a}ticas, Universidad Nacional Aut\'{o}noma de M\'{e}xico}, Javier Bracho\thanks{Instituto de Matem\'{a}ticas, Universidad Nacional Aut\'{o}noma de M\'{e}xico}, Natalia Garc\'{i}a-Col\'{i}n\thanks{Instituto de Matem\'{a}ticas, Universidad Nacional Aut\'{o}noma de M\'{e}xico}, Isabel Hubard\thanks{Instituto de Matem\'{a}ticas, Universidad Nacional Aut\'{o}noma de M\'{e}xico}}

%\address[UNAM]{Instituto de Matem\'{a}ticas, Universidad Nacional Aut\'{o}noma de M\'{e}xico, \'{A}rea de la Investigaci\'{o}n Cient\'{i}fica, Circuito Exterior, Ciudad Universitaria, Coyoac\'{a}n, 04510, M\'{e}xico D.F., M\'{E}XICO}

\maketitle

\begin{abstract}

The \emph{intersection matrix} of a finite simplicial complex has as each of its entries the rank of the intersection of its respective simplices. We prove that such matrix defines the triangulation of a closed connected surface up to isomorphism.

\end{abstract}

\section{Introduction \& Motivation}

Within the theory of convex polytopes, the study of the combinatorial equivalence of $k$-skeleta of pairs polytopes which are not equivalent themselves has been of interest, this phenomena is referred to in the literature as ambiguity \cite{grunbaum1967convex}.

It is well known that for $k \geq \lfloor \frac{d}{2} \rfloor$ the $k$-skeleton of a convex polytope is not dimensionally ambiguous, this is, it defines the entire structure of its underlying $d$-polytope. However for $k < \lfloor \frac{d}{2} \rfloor$ the question is much more intricate.

One of the most interesting results in this direction is the solution to Perle's conjecture by P.Blind and R.Mani \cite{Blind1987} and, separately, by G. Kalai \cite{kalai1988simple} which states that the 1-skeleta of convex simple $d$-polytopes define their entire combinatorial structure. Or, on its dual version, that the dual graph (facet adjacency graph) of a convex simplicial d-polytope determines its entire combinatorial structure. (Also see \cite{ziegler1995lectures})

In its simplest version this theorem states that a triangulation of a convex polyhedra is defined up to isomorphism by its dual graph. We were motivated by this assertion to reveal possible generalizations of that theorem.

\section{Contribution}

A \emph{triangulated surface} is a simplicial complex whose underlying topological space is a connected $2$-manifold (without boundary). It is known that the dual graph of a triangulated space does not define it. In fact, there exist different triangulated surfaces having the same dual graph (see Goldberg snarks on \cite{mohar2004polyhedral}, pages 16-18). However, as building appropiate examples of the latter is not simple, this strongly suggests that more information is needed in order to reconstruct triangulations of surfaces (up to isomorphism). Here we exhibit exactly what additional information is needed, and argue in the conclusions that our hypothesis are optimal for these surfaces.

For a triangulated surface $S$ we will denote $V, E$ and $T$ as its sets of vertices, edges and triangles, respectively. We will say that two triangulations $S$ and $S'$ have the same intersection matrix if there is a bijective map, $f: T \rightarrow T'$, from the triangles of $S$ into the triangles of $S'$  such that for any two triangles $t_1, t_2 \in T$, the equality $| t_1\cap t_2|=| f(t_1)\cap f(t_2)|$ holds. Using this terminology we can now state our main result as follows:

\begin{theorem} \label{thm:main} Two triangulated surfaces which have the same intersection matrix are isomorphic.
\end{theorem}

This result will follow as a corollary of a more precise statement for which we need to introduce a few further definitions.

\section{Preliminaries}

\begin{figure}[t]
\begin{center}
\includegraphics[scale=0.55]{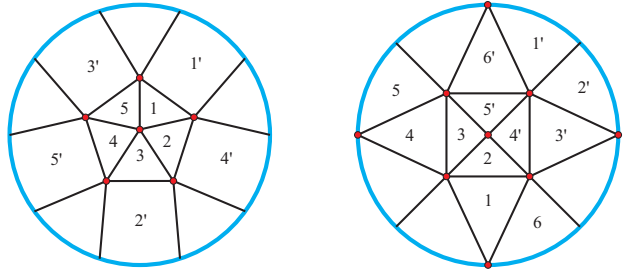} 
\label{fig:PT}
\caption{The only two triangulated surface which have intersection preserving maps not extendable to isomorphisms.}
\end{center}
\end{figure}

Given two triangulations of a surface  $S$ and $S'$, for convenience, we will say that a map $f: T \rightarrow T'$, between the triangles of $S$ and $S'$ which preserves their intersection matrix is an \emph{intersection preserving mapping}, additionally we will say that $f$ \emph{extends to an isomorphism} if there is a simplicial map $g:V \rightarrow V'$ which induces $f$. Not every intersection preserving map can be extended to an isomorphism. The triangulations of the projective plane shown in Figure \ref{fig:PT} are examples of such an occurrence.

The triangulation on the left hand side in the figure is the \emph{half-icosahedron} which we will denote $\mathcal{I}/2$, while the triangulation on the right hand side in the figure is a triangulation of the \emph{half-cube}, where each square facet is subdivided into four triangles by adding one additional vertex at its centre. This triangulation will be denoted as $\mathcal{TC}/2$. It can be easily checked that the intersection preserving map given by the labels of the triangles of each one of the triangulations into itself ($i \leftrightarrow i'$) is not extendable to an isomorphism.

In essence Theorem \ref{thm:main} holds because the half-icosahedron and triangulated half-cube are the only two examples of triangulations of closed surfaces which have intersection preserving mappings that cannot be extended to isomorphisms. More precisely:

\begin{theorem} \label{thm:main2} Let $S$ and $S'$ be two triangulated connected closed surfaces with no boundary and let $f: T \rightarrow T'$ be an intersection preserving map between their respective sets of triangles which does not extend to an isomorphism; then either $S \simeq S' \simeq \mathcal{I}/2$ or $S \simeq S' \simeq \mathcal{TC}/2$.
\end{theorem}

One of the characteristics of triangulations of a closed surface is that the neighbourhood of every vertex is a disk. Furthermore, the triangles incident to any vertex of such surface form the simplest of triangulations of a disk, namely an $n$-gon whose vertices are all linked by an edge to a central vertex in the centre of the $n$-gon. We start off by analysing the intersection patterns of such a structure, and find the two other objects which can have an intersection matrix equal to that of a triangulated disk.
 
\begin{definition} An n-shell is the abstract triangulation $\Delta$ such that it set of triangles $T_\Delta=\{t_0, t_1,\ldots, t_{n-1}\}$ satisfies $|t_i \cap t_{i+1\ mod n}|=2$ and $|t_i \cap t_{j}|=1$ for $|i-j|\geq2$ with $i,j \in \{0,\ldots, n-1\}$.
\end{definition}

\begin{lemma} \label{lem:disk} The vertex sets of the triangles in an $n$-shell, with $n\geq 3$, $T_\Delta=\{t_0, t_1,\ldots, t_{n-1}\}$ can be equivalent to the vertex sets of a triangulation of one of the following three objects;
\begin{enumerate}

\item a triangulated disk (i.e. an $n$-shell),\\
$t_{i}=\{a_{i \mod n}, a_{{i+1} \mod n}, x\}$ for all $0\leq i \leq n$, for any $n$;

\item a triangulation of a M\"{o}bius band with five triangles, \\
$t_{0}= \{a_{0}, a_{2}, a_{1}\}$, 
$t_{1}= \{a_{1}, a_{3}, a_{2}\}$,
$t_{2}= \{a_{2}, a_{4}, a_{3}\}$,
$t_{3}= \{a_{3}, a_{0}, a_{4}\}$, and
$t_{4}= \{a_{4}, a_{1}, a_{0}\}$; or

\item a triangulation of a M\"{o}bius band with six triangles, \\
$t_{0}=\{a_{0}, a_{1}, a_{2}\}$, $t_{1}= \{a_{1}, a_{2}, a_{4}\}$, 
$t_{2}= \{a_{2}, a_{3}, a_{4}\}$,
$t_{3}= \{a_{3}, a_{0}, a_{4}\}$,
$t_{4}= \{a_{0}, a_{5}, a_{4}\}$, and
$t_{5}= \{a_{5}, a_{2}, a_{0}\}$.

\end{enumerate}
\end{lemma}

The proof of this lemma consists of several parts and follows largely by a detailed analysis of the combinatorial structure of $n$-shells of triangles. We include it in detail in Appendix \ref{ap:lem}.

%%%%%%%%%%%%%%%%%%%%%%%%%%%%%%%%%%%%
%%%%%%%%%%%%%%%%%%%%%%%%%%%%%%%%%%%%%%
\section{Proof Theorem \ref{thm:main2}}

We now proceed to proving Theorem 2, using exhaustively the local and global implications of Lemma \ref{lem:disk}.
\begin{proof}

For each vertex $v \in V$, the set of vertices of the triangulation $S$, let $\Delta$ be the $n_{v}$-shell around $v$, by hypothesis union of the triangles in $\Delta$ is necessarily a disk. 

Note that if for all $v \in V$, the triangles $f(T_\Delta)$ in $S'$ also form a disk, then the mapping $h: V \rightarrow V'$, from the vertices of $S$ into the vertices of $S'$ such that $h(v)= \bigcap_{t \in T_\Delta} f(t)$ is a simplicial mapping which induces $f$. Hence we can assume that there is a vertex $v \in V$ such that the union of the triangles in $f(\Delta)$ is not a disk.

\emph{Case 1.} Suppose the union of the triangles $f(T_\Delta)$ is the $5$-triangulation of the M\"{o}bius band described in Lemma \ref{lem:disk}.

Let $T_\Delta= \{ t_{0}, t_{1}, t_{2}, t_{3}, t_{4}\}$, where $t_{i}=\{a_{i \mod 5}, a_{i+1 \mod 5}, x\}$ and  $f(T_\Delta)=\{ t_{0}', t_{1}', t_{2}', t_{3}', t_{4}'\}$ where 
$t_{0}'= \{a_{0}', a_{2}', a_{1}'\}$,
$t_{1}'= \{a_{1}', a_{3}', a_{2}'\}$,
$t_{2}'= \{a_{2}', a_{4}', a_{3}'\}$,
$t_{3}'= \{a_{3}', a_{0}', a_{4}'\}$, and
$t_{4}'= \{a_{4}', a_{1}', a_{0}'\}$.

Given that $S'$ is also a closed connected surface, then each of the simplices $ t_{0}', t_{1}', t_{2}', t_{3}', t_{4}'$ has got a triangle adjacent to its remaining free edge. Let $r'_{i}$ be the simplices such that $|r'_{i}\cap t'_{i}|=2$, then 
$r_{0}'= \{a_{0}', a_{2}', x_{0}'\}$,
$r_{1}'= \{a_{1}', a_{3}', x_{1}'\}$,
$r_{2}'= \{a_{2}', a_{4}', x_{2}'\}$,
$r_{3}'= \{a_{3}', a_{0}', x_{3}'\}$, and
$r_{4}'= \{a_{4}', a_{1}', x_{4}'\}$. \\ This  is
$r_{i}'= \{a_{i\mod 5}', a_{i+2 \mod 5}', x_{i}'\}$. It follows that, $|r_{i}' \cap t_{j}'|\geq 1$ for all $i\neq j$. 

Note that the interior of each of the edges $\{a_{i\mod 5}', a_{i+1 \mod 5}'\}$  is in the interior of the M\"{o}bius band, thus this edges cannot be repeated in any further simplex in the complex. This implies that $x_{i} \not \in \{a_{0}', a_{1}', a_{2}', a_{3}', a_{4}'\}$, because, if this was the case, at least one of the edges $\{a_{i\mod 5}', a_{i+1 \mod 5}'\}$ would be a subset of $r'_{i}$. Then, $|r_{i}' \cap t_{j}'|= 1$.

Let $r_{i}=f^{-1}(r'_{i})$, then $|r_{i} \cap t_{i}|=2$ and $|r_{i} \cap t_{j}|= 1$ for all $i\neq j$. As $t_{i}=\{a_{i \mod 5}, a_{i+1 \mod 5}, x\}$ then $r_{i}=\{a_{i \mod 5}, a_{i+1 \mod 5}, x_{i \mod 5}\}$. Here $|r_{i\mod 5} \cap t_{i+1\mod 5}|\geq 1$ and $|r_{i\mod 5} \cap t_{i-1\mod 5}|\geq 1$ trivially, hence \\ $a_{i-1 \mod 5}, a_{i+2 \mod 5} \not \in r_{i}$. However, for $|r_{i\mod 5} \cap t_{i+2\mod 5}|=1$ and $|r_{i\mod 5} \cap t_{i-2\mod 5}|=1$ to be accomplished, necessarily $x_{i}=a_{i+3\mod 5}=a_{i-2\mod 5}$. That is, $r_{i}=\{a_{i \mod 5}, a_{i+1 \mod 5}, a_{i-2 \mod 5}\}$, hence the simplicial complex asociated to $\bigcup_{i=0}^{4}r_{i}$ is a $5$-triangulation of a M\'{o}bius band, where $|r_{i} \cap r_{i+2 \mod 5}|=|r_{i} \cap r_{{i-2} \mod 5}|=2$, thus $T = \bigcup_{i=0}^{4}r_{i} \cup \bigcup_{i=0}^{4}t_{i}$ and the triangulations associated to $\bigcup_{i=0}^{4}r_{i}$ and $\bigcup_{i=0}^{4}t_{i}$ are a M\"{o}bius band and a disk, respectively, then it is straightforward to check that $S$ is the triangulation $\mathcal{I}/2.$

The latter also implies that $|r'_{i} \cap r'_{i+2 \mod 5}|=|r'_{i} \cap r'_{i-2 \mod 5}|=2$ then $v'=v'_{i}$ for all $i=1,\ldots 4$, so that $T' =\bigcup_{i=0}^{4}r'_{i} \cup \bigcup_{i=0}^{4}t'_{i}$, hence $S'$ is also the triangulation $\mathcal{I}/2$.

\emph{Case 2.} Suppose the union of the triangles $f(T_\Delta)$ is the $6$-triangulation of the M\"{o}bius band described in Lemma \ref{lem:disk}.

Let $T_\Delta= \{ t_{0}, t_{1}, t_{2}, t_{3}, t_{4}, t_{5}\}$, where $t_{i}=\{a_{i \mod 6}, a_{i+1 \mod 6}, x\}$ and  $f(T_\Delta)=\{ t_{0}', t_{1}', t_{2}', t_{3}', t_{4}', t_{5}'\}$ where 
$t'_{0}=\{a'_{0}, a'_{1}, a'_{2}\}$,
$t'_{1}= \{a'_{1}, a'_{2}, a'_{4}\}$,
$t'_{2}= \{a'_{2}, a'_{3}, a'_{4}\}$,
$t'_{3}= \{a'_{3}, a'_{0}, a'_{4}\}$,
$t'_{4}= \{a'_{0}, a'_{5}, a'_{4}\}$, and
$t'_{5}= \{a'_{5}, a'_{2}, a'_{0}\}$.

As $S'$ is a connected closed surface, then each of the simplices $\{ t_{0}', t_{1}', t_{2}', t_{3}', t_{4}', t_{5}'\}$ has got a triangle adjacent to its remaining free edge. Let $r'_{i}$ be the simplices such that $|r'_{i}\cap t'_{i}|=2$, then 
$r_{0}'= \{a_{0}', a_{1}', x_{0}'\}$,
$r_{1}'= \{a_{1}', a_{4}', x_{1}'\}$,
$r_{2}'= \{a_{2}', a_{3}', x_{2}'\}$,
$r_{3}'= \{a_{3}', a_{0}', x_{3}'\}$,
$r_{4}'= \{a_{4}', a_{5}', x_{4}'\}$, and
$r_{5}'= \{a_{2}', a_{5}', x_{5}'\}$.

Here it follows that, $|r_{i}' \cap t_{j}'|\geq 1$ for all $i\neq j$, except for the pairs $i=0 \text{ and } j=2$, $i=1\text{ and }j=5$, $i=2\text{ and }j=4$, $i=3\text{ and }j=1$, $i=4\text{ and }j=0$ and $i=5\text{ and }j=3$; for these exceptions the intersection might be empty.

The latter implies that, if $r_{i}=f^{-1}(r'_{i})$, then $|r_{i} \cap t_{i}|=2$ and $|r_{i} \cap t_{j}|\geq 1$ for all $i\neq j$, except for the pairs $i=0 \text{ and } j=2$, $i=1\text{ and }j=5$, $i=2\text{ and }j=4$, $i=3\text{ and }j=1$, $i=4\text{ and }j=0$, and $i=5\text{ and }j=3$; for these exceptions the intersection might be empty.

As $t_{i}=\{a_{i \mod 6}, a_{i+1 \mod 6}, x\}$ then the vertex sets of the $r_{i}$'s are $r_{i}=\{a_{i \mod 6}, a_{i+1 \mod 6}, x_{i \mod 6}\}$. Note that $x_{0}' \not \in \{a'_{0}, a'_{1}, a'_{2}, a'_{4}, a'_{5}\}$ as the edges $\{a'_{0}, a_{2}\}$, $\{a'_{0}, a_{4}\}$, $\{a'_{0}, a_{5}\}$, $\{a'_{0}, a_{4}\}$ are edges whose interior is in the interior of the M\"{o}bius band. Thus we might have $v'_{0}=a'_{3}$, however if that was the case
$|r_{0}' \cap r_{3}'|= 2$, and $|r_{0} \cap r_{3}|= 2$, but this is not possible. Then necessarily  $x_{0}' \not \in \{a'_{0}, a'_{1}, a'_{2}, a'_{3}, a'_{4}, a'_{5}\}$. 

Using an argument analogous to the one in the previous case, we deduce that for each $i$, $x_{i}' \not \in \{a'_{0}, a'_{1}, a'_{2}, a'_{3}, a'_{4}, a'_{5}\}$; so that $|r_{i}' \cap t_{j}'|= 1$ for all $i\neq j$, except for the pairs $i=0\text{ and }j=2$, $i=1\text{ and }j=5$, $i=2\text{ and }j=4$, $i=3\text{ and }j=1$, $i=4\text{ and }j=0$, and $i=5\text{ and }j=3$, for which the intersection is empty. 

The above implies $|r_{i} \cap t_{i}|=2$ and $|r_{i} \cap t_{j}|= 1$ for all $i\neq j$, except for the pairs $i=0\text{ and } j=2$, $i=1\text{ and } j=5$, $i=2\text{ and } j=4$, $i=3\text{ and } j=1$, $i=4\text{ and } j=0$, and $i=5\text{ and } j=3$, for which the intersection is empty. Hence, in order to accomplish the intersection dimensions indicated by the map necessarily, $x_{0}=x_{1}= a_{4}$, $x_{2}=x_{3}= a_{0}$, $x_{4}=x_{5}= a_{2}$, thus;
$r_{0}=\{a_{0}, a_{1}, a_{4}\}$,
$r_{1}=\{a_{1}, a_{2}, a_{4}\}$,
$r_{2}=\{a_{2}, a_{3}, a_{0}\}$,
$r_{3}=\{a_{3}, a_{4}, a_{0}\}$,
$r_{4}=\{a_{4}, a_{5}, a_{2}\}$, and
$r_{5}=\{a_{0}, a_{5}, a_{2}\}$.

Therefore, the triangulation associated to $\bigcup_{i=0}^{5} r_{i}$ is a $6$-triangulation of a M\"{o}bius band and $T= \bigcup_{i=0}^{5} t_{i} \cup \bigcup_{i=0}^{5} r_{i}$. It is now straightforward to check that $S$ is equal to the triangulation $\mathcal{TC}/2$.

The implication for $S'$ is that $|r'_{0} \cap r'_{1}|=2,\; |r'_{1} \cap r'_{4}|=2, \; |r'_{4} \cap r'_{5}|=2, \; |r'_{5} \cap r'_{2}|=	2, \; |r'_{2} \cap r'_{3}|=2, \; |r'_{3} \cap r'_{0}|=2$, which in turn implies $v'=v'_{i}$ for all $i \in \{0,\ldots 5\}$ and, further, $T'= \bigcup_{i=0}^{5} t'_{i} \cup \bigcup_{i=0}^{5} r'_{i}$, so that $S'$ is also equal to the triangulation $\mathcal{TC}/2$.

\end{proof}
%%%%%%%%%%%%%%%%%%%%%%%%%%%%%%%%%%%%%%%%%%%%%%%%%%

\section {Proof Lemma \ref{lem:disk}} \label{ap:lem}

The proof  will proceed by induction on the number of triangles on the $n$-shell. Before we proceed with the proof we need to define a special set of vertices present in every $n$-shell.

\begin{definition} The structural vertex list of an $n$-shell, $\Delta$ such that $T_\Delta=\{t_0, t_1,\ldots, t_{n-1}\}$, is the vertex set $\{a_{0}, a_{1}, \ldots, a_{n-1}, b_{0}, b_{1}, \ldots, b_{n-1}\} \subset \;_{n}\Pi^{2}_{0}$, where $a_{i}=t_{i}\setminus t_{i+1 \mod n }$ and $b_{i}=t_{i}\setminus t_{i-1 \mod n}$ for $0\leq i\leq n-1$.

If the $n$-shell is open (i.e. $|t_1 \cap t_{n-1}|=1$) then the structural vertex list is only defined for $0\leq i\leq n-2$ and $1\leq i\leq n-1$.
\end{definition}

Note that there might be some vertex repetition in the structural vertex list of any given n-shell; however, by definition, we can be certain that 
\begin{equation} \label{eqn:begg}
a_{i}\neq b_{i} \text{ for all } 0 \leq i \leq n-1 \text{ fixed,}
\end{equation} 
\begin{equation}
a_{i}\neq a_{i+1 \mod n} \text{\; and\; } b_{i}\neq b_{i+1 \mod n} \text{ for all } 0 \leq i \leq n-1.
\end{equation}
Moreover, by definition, neither $a_{i}$ nor $b_{i+1} \in \; t_{i}\cap t_{i+1}$ for any $0\leq i \leq n-1$ thus
\begin{equation} \label{eqn:end}
a_{i}\neq b_{i+1} \text{ for all } 1\leq  i < n-1.
\end{equation} 

For any $i$, observe that $t_{i-1 \mod n}=\{a_{i-1 \mod n}, a_{i}, x_{i}\}$, $t_{i}=\{a_{i}, b_{i}, x_{i}\}$ and $t_{i+1\ mod n}=\{b_{i}, b_{i+1 \mod n}, x_{i}\}$, where $x_{i}= t_{i-1}\cap t_{i} \cap t_{i+1}$.

In order to prove Lemma 1 we will study the repetition patterns of the structural vertex list and show that they characterize the different combinatorial types of $n$-shells, through the use of the following claims:

\begin{claim}There is no vertex repetition in the structural vertex list of any open 3-shell, $\{a_{i-1},a_{i},b_{i}, b_{i+1}\}$.
\end{claim}
By equations (1), (2) and (3) it only remains to be proven that $a_{i-1 \mod n} \neq b_{i+1 \mod n}$. If that was the case then  $(t_{i-1 \mod n})_{0}=\{a_{i-1 \mod n}, a_{i}, x_{i}\}= (t_{i+1\ mod n})_{0}=\{b_{i}, b_{i+1 \mod n}, x_{i}\}$, forcing $n=3$.

\begin{claim}The vertices $a_{i}\neq a_{j}$  and $b_{i}\neq b_{j}$, for all $i, j \in \{0,\ldots,n-1\}$.
\end{claim}
First observe that the statement follows directly for $n$-shells if we prove it for open $n$-shells. The proof will follow by induction on the distance $|i-j|$ between the indices we are comparing.

Equation (2) proves the statement for $|i-j|=1$, while Claim 1 proves the statement for $|i-j|=2$. We assume the statement holds for $|i-j|\leq k>2$. Let $|i-j|=k+1$.

First suppose $a_i=a_j$. Considering that $ t_{i-1}=\{a_{i-1}, a_{i}, x_{i}\}$, $t_{i}=\{a_{i}, b_{i}, x_{i}\}$; and $t_{j-1}=\{a_{j-1}, a_{j}, x_{j}\}$, $t_{j}=\{a_{j}, b_{j}, x_{j}\}$, then$|t_{i-1}\cap t_{j-1}|= 1$, $|t_{i-1}\cap t_{j}|=1$, $|t_{i}\cap t_{j-1}|=1$ and $|t_{i}\cap t_{j}|=1$ already happen with just one repetition, $\{a_{i-1}, x_{i}, b_{i}, a_{j-1}, x_{j}, b_{j}, a_i=a_j\}$ are necessarily all different among them.

Given that $t_{j+1}\cap t_{i-1} =\{a_{i-1}, a_{i}, x_{i}\}\cap \{b_{j}, b_{j+1}, x_{j}\}$ and $t_{j+1}\cap t_{i}=\{b_{j}, b_{j+1}, x_{j}\}\cap \{a_{i}, b_{i}, x_{i}\}$, and that the intersection of the two pairs of sets above is of cardinality one, necesarily $b_{j+1}=a_i$ or $b_{j+1}=x_i$. Also, as $t_{i+1}\cap t_{j} =\{a_{j}, b_{j}, x_{j}\}\cap \{b_{i}, b_{i+1}, x_{i}\}$  have intersection of cardinality one too, so either $b_{i+1}=a_j$ or $b_{i+1}=x_j$. However the repetitions in this paragraph and last paragraph can't both happen simultaneously, preserving all the dimensions of the intersections.

The proof for $b_{i}\neq b_{j}$ is analogous.

\begin{claim}  The structural vertex list of every four consecutive triangles, $t_i, t_{i+1},$ $t_{i+2}, t_{i+3}$, in an $n$-shell with $n \geq 5$ has exactly one repetition which is either $b_{i+1 }=a_{i+2}$ or $a_{i }=b_{i+3}$.
\end{claim}

By Claim 1 the lists $a_{i}, a_{i+1}, b_{i+1}, b_{i+2}$ and $a_{i+1}, a_{i+2},  b_{i+2}, b_{i+3}$ considered separately have no repetition among them and by Claim 2, all a's are different to the b's; thus, the only possible repetitions are $a_{i}=b_{i+3}$, or $b_{i+1}=a_{i+2}$, and they can't both occur simultaneously whilst maintaining the dimension of the intersections of the triangles.

\begin{claim} \label{claim:four} The vertex sets of four consecutive triangles $t_i, t_{i+1}, t_{i+2}, t_{i+3}$ in an $n$-shell are completely determined by the repetition of their structural vertex list, according to the following rules;
\begin{enumerate}[(a)]
\item if  $b_{i+1 }=a_{i+2}$ then $t_{i}=\{a_{i}, a_{i+1}, x\}$, 
$t_{i+1}=\{a_{i+1}, b_{i+1}=a_{i+2}, x\}$, 
$t_{i+2}=\{a_{i+2}=b_{i+1}, b_{i+2}, x\}$ and 
$t_{i+3}=\{b_{i+2}, b_{i+3}, x\}$;
\item if  $a_{i}=b_{i+3}$ then $t_i=\{a_i=b_{i+3}, a_{i+1}, a_{i+2}\}$,
$t_{i+1}=\{a_{i+1}, b_{i+1}, a_{i+2}\}$
$t_{i+2}=\{b_{i+1}, b_{i+2}, a_{i+2}\}$ and 
$t_{i+3}=\{b_{i+2}, b_{i+3}=a_i, b_{i+1}\}$.
\end{enumerate}
\end{claim}

In the first case (a) as, $t_i=\{a_i, b_{i}, x_{i}\}$ and $b_{i+1 }=a_{i+2}$, then 
$t_{i+1}=\{a_{i+1}, b_{i+1}=a_{i+2}, x_{i+1}\}=\{a_{i+2}=b_{i+1}, a_{i+1}, x_{i+2}\}$ implies $x_{i+1}=x_{i+2}$ and the statement follows.

In the second case (b), if $a_{i}=b_{i+3}$ then we then have $t_i=\{a_i=b_{i+3}, a_{i+1}, x_{i+1}\}$
$t_{i+1}=\{a_{i+1}, b_{i+1}, x_{i+1}\}=\{a_i=b_{i+3}, a_{i+1}, x_{i+2}\}$
$t_{i+2}=\{b_{i+1}, b_{i+2}, x_{i+1}\}=\{a_{i+2}, b_{i+2}, x_{i+2}\}$ and 
$t_{i+3}=\{b_{i+2}, b_{i+3}=a_i, x_{i+2}\}$. Then as $|t_i \cap t_{i+3}|=1$ all of  $a_{i+1}, x_{i+1}, b_{i+2}, x_{i+2}$ are different. Then the two different expressions of $t_{i+1}$ and $t_{i+2}$ imply $b_{i+1}=x_{i+2}$ and $a_{i+2}=x_{i+1}$, respectively and the statement follows.

\begin{figure}[h]
\begin{center}
\includegraphics[scale=0.8]{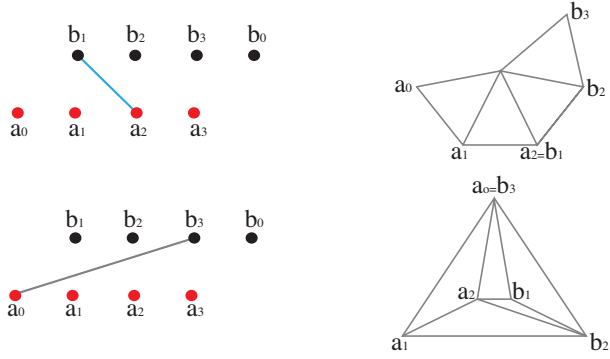} \label{fig:openshell}
\caption{The two possible triangulations associated to an open $4$-shell.}
\end{center}
\end{figure}

Figure \ref{fig:openshell} depicts the two possible repetitions in the structural vertex list of four triangles, as per Claim \ref{claim:four}; and the resulting sets of vertices of the triangles.

 \begin{claim}The two occurring repetitions in the structural vertex list of the first five consecutive triangles of an open $n$-shell ($n\geq 5$) determine the repetition pattern of the structural vertex list of the whole $n$-shell, furthermore 
the repetition pattern is one of the following:
\begin{enumerate}[(a)]
\item $b_{1}=a_{2}$ and $b_{2}=a_{3}$, and  $b_{i}=a_{i+1}$ for all $i \in \{1, \ldots, n-3\}$;
\item $b_{1}=a_{2}$ and $a_{1}=b_{4}$, then $n \leq 6$ and when $n=6$, $b_{3}=a_{4}$;
\item $a_{0}=b_{3}$ and $b_{2}=a_{3}$, then $n \leq 6$ and when $n=6$, $a_{2}=b_{5}$;
\item $a_{0}=b_{3}$ and $a_{1}=b_{4}$, then $n \leq 4$.
\end{enumerate}
\end{claim}

First observe that no pair of repetitions, as allowed by Claim 3, among the groups $b_{1}=a_{2}$, $a_{2}=b_{5}$ and  $b_{3}=a_{4}$ , $a_{0}=b_{3}$, might ever happen simultaneously. For instance, if $b_{1}=a_{2}=b_{5}$ then two b's would be equal in contradiction Claim 2. 

Then we might have that two consecutive repetitions in the list can only possible be those in the statement of the claim.

(a) For $n=5$, the statement holds trivially. The proof will follow by induction on the number of triangles n. Assume the statement holds for $n\leq k$, and let $5 < n= k+1$. By the induction hypothesis $b_{i}=a_{i+1}$ for all $i \in \{1, \ldots, k-3\}$ and by Claim 4, we can deduce $b_{k-2}=a_{k-1}$, thus the statement follows.

(b) For $n=5$, the statement holds trivially. Consider the first five simplices of the n-shell. Then by Claim 4
$t_{0}=\{a_{0}, a_{1}=b_4, x\}$,\\
$t_1=\{a_1=b_{4}, a_{2}=b_1, x\}=\{a_1=b_{4}, a_{2}=b_1, a_{3}\}$,\\
$t_{2}=\{a_{2}=b_1, b_{2}, x\}=\{a_{2}=b_1, b_{2}, a_{3}\}$,\\
$t_{3}=\{b_{2}, b_{3}, x\}=\{b_{2}, b_{3}, a_{3}\}$, and \\
$t_{4}=\{b_{3}, b_{4}=a_1, b_{2}\}$, so that $x=a_3$ then all intersection sizes are fulfilled, hence no other repetition can happen. 

Consider now $t_2, t_3, t_4, t_5$ then $t_5=\{a_1=b_4, b_5, x_4\}$. Here $x_4$ is either equal to $a_4=b_3$ or $b_2$. In the first case, $b_5$ cannot be made equal to any other vertex in order to complete the dimension of the intersection of $t_5$ with $t_2$ whilst fulfilling all the other dimensions of intersection.
Hence $x_4=a_4=b_3$, all the dimensions of the intersections are fulfilled, and no other repetition can occur.

Finally, consider $t_6$ and its vertex set  $t_6=\{b_5, b_6, x_5\}$,  then $x_5$ can be either $b_4=a_1$ or $b_2$. However in both cases it is impossible to complete the dimensions of the pairwise intersections with the other triangles. Thus, we have proven that in this case $n\leq 6$, and the statement follows.

(c) Follows in an analogous manner to (2).

(d) If $n=5$, the statement follows trivially, and by Claim 4,
$t_0=\{a_0=b_3, a_1, a_2\}$, $t_1=\{ a_1, a_2, a_3\}$, $t_2=\{a_2, b_2, a_3\}$,  $t_3=\{a_0=b_3, b_2, a_3\}$, $t_4=\{a_0=b_3, b_2, b_4\}$.  However the dimension of the intersection of $t_4$ with $t_0$ and $t_1$ has not been fulfilled. The only possible manner for the respective two intersections to be fulfilled without altering the rest of the intersections is for $a_1=b_4$. But in this case $dim(t_4 \cap t_0)=1$. Thus, we have proven that in this case $n\leq 4$, and the statement follows.

Using the previous claims the statement of the lemma follows easily. If $n\leq 4$, it can be easily verified using Claims 1 and 4.

For $n \geq 5$,  clearly,  case 1 in the statement of the lemma corresponds to case (a) in Claim 5, for any $n$.

It is also easy to check that there are no (closed) $5$-shells that correspond to the repetitions of the structural vertex list in cases (b) and (c) in Claim 5. However there is a unique (closed) $6$-shell corresponding to the repetition types of cases (b) and (c), this has the vertex set of case 3 in the statement of the lemma.

Finally, within the proof of case (d) in Claim 5 it was shown that there is only one (closed) $5$-shell with this type of repetiton on its vertex list. This $5$-shell has vertex sets as per case 2 in the statement of the lemma. \qed

\begin{figure}[h]
\begin{center}
\includegraphics[scale=0.8]{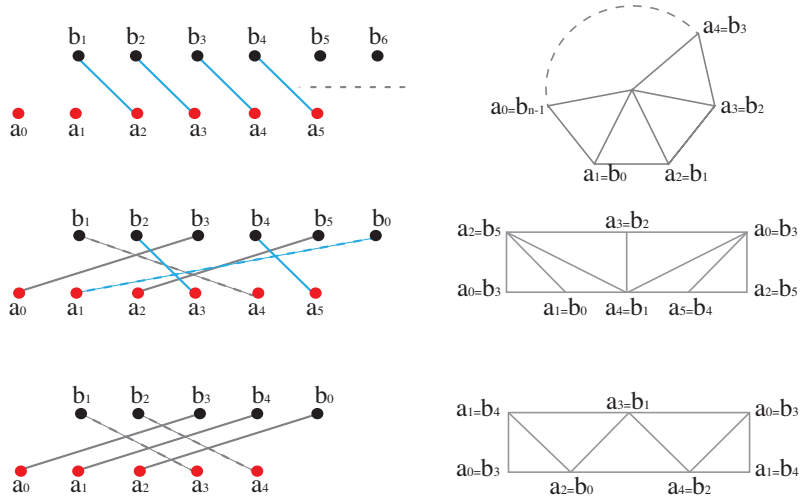} \label{fig:closedshell}
\caption{The three possible traingulations associated to an $n$-shell. }
\end{center}
\end{figure}

Figure \ref{fig:closedshell} depicts the three possible repetitions of the structural vertex list of a closed $n$-shell and the corresponding resulting sets of vertices for a triangulation with shuch repetitions; these correspond to the three triangulations listed in Lemma \ref{lem:disk}.

%%%%%%%%%%%%%%%%%%%%%%%%%%%%%%%%%%%%%%%%%%%%%%%%%%%%%%%
\section{Conclusion}

The proof of Theorem \ref{thm:main} follows directly from Theorem \ref{thm:main2}. 

As mentioned in the introduction, it is known that the dual graph of a triangulated space does not define it and in fact, there exist different triangulated closed connected surfaces with no boundary having the same dual graph \cite{mohar2004polyhedral}. The reader might wonder if the intersection matrix of a triangulation characterizes a connected surface with a boundary. The answer to this question is negative and some examples of this are provided by the objects in Lemma \ref{lem:disk}. 

However some of the surfaces in the lemma are non orientable, thus one may further ask if the intersection matrix of a connected orientable surface defines it. Here, we conjecture that the answer to this question is possitive, as all examples of connected orientable surfaces having the same intersection matrix which we have found, are isomorphic, but, just as for the triangulations of the projective plane of Theorem \ref{thm:main2}, they do not have intersection preserving maps which extend to isomorphisms between them.

Further, one might wonder on the generalizations of Theorem \ref{thm:main} to simplicial complexes of any rank. We have strong evidence to suggest that the statement of Theorem \ref{thm:main} can indeed by extended to any rank, provided that we allow for some complexes which have no intersection preserving maps that extend to isomorphisms between them, but are isomoprhic nonetheless.

%%%%%%%%%%%%%%%%%%%%%%%%%%%%%%%%%%%%%%%%%%%%%%%%%%

%%%%%%%%%%%%%
%BIBLIOGRAPHY
\bibliographystyle{plain}
\bibliography{simplicialpuzzles}
\end{document}